# Simultaneous Detection and Estimation, False Alarm Prediction for a Continuous Family of Signals in Gaussian Noise

D. Michael Milder, Robert G. Lindgren, and Morris M. Berman

*Abstract*—New problems arise when the standard theory of joint detection and estimation is applied to a set of signals drawn from a continuous family; decision thresholds must be determined as a function of the continuous parameter $x$ characterizing the signals, and false alarms occur, not with a discrete probability, but with a density in $x$. A Bayes decision structure over the domain of signal parameters yields a state estimate of the signal parameter $x$ as an integral part of a signal declaration. The decision criterion is converted to a form in which detection and false alarm densities appear and from which is derived a relation between them for all $x$. The limiting case of additive Gaussian noise and a high detection threshold allows a simplified decision criterion and a state estimate of signal location in $x$ that approaches the Cramer-Rao bound. Also in this limit, an analytic form for the false alarm probability density over $x$, a quantity not readily obtained in general, is evaluated here through its relation to the detection probability. The false alarm density expression and state accuracy prediction are tested through Monte Carlo simulations, and the comparison demonstrates excellent agreement.

## I. INTRODUCTION

Detection of known signals in noise is seldom unaccompanied by estimation. The unknown vector of parameters $x$ characterizing the family of signals commonly occurs in its simplest form as a time of arrival or a location in an image, and the operation of signal detection requires, at least implicitly, a search over the parameter space. Middleton and Esposito [1] pioneered the unification of detection and estimation as a problem in Bayesian decision theory, examining several scenarios of sequential decisions with both independent and coupled cost functions for the two operations. Subsequent papers have elaborated on this fundamental approach, for Gaussian noise [2], for bounds on detection and estimation performance [3], and for the advantage of joint, rather than sequential optimization [4]. The parameter estimation decisions in these papers, however, are discrete. Moustakides et al. extend the analysis to continuous signal parameters and derive "decision structures" to minimize some variant of the Bayes risk (misestimation cost plus missed-signal cost in [5], misestimation cost conditional on a signal-present decision in [6]) subject to a false alarm constraint (constraint on Bayes cost under a no-signal hypothesis in [7]). In this paper we also approach the optimization of simultaneous detection and estimation through the minimization of the Bayes risk, but we avoid the complication caused by a false alarm constraint by assigning Bayesian costs to all decisions. A decision structure is then readily obtained that minimizes the Bayes risk over all measurement possibilities. This setup represents the starting point for our analysis, which focuses on the determination of signal detection probabilities, state estimation errors, and false alarm rates as functions of the signal parameter $x$. Variation of the false alarm cost, which we shall use to derive an important relation between signal detection and false alarm probabilities as functions of the signal parameter $x$, is the effective equivalent to variation of a false-alarm constraint.

We begin by specifying an optimization criterion based on the Bayes risk, and from it we construct a decision statistic $T$ in which the likelihood ratio for signal versus noise, formed for a possible signal at any $x$, has a central role. When this methodology is applied to the problem of a continuous family of signals in additive Gaussian noise, and the estimate of $x$ is allowed to be continuous, some novel and useful properties emerge. The decision test at high thresholds (high false-alarm costs) can be cast in a form that is greatly simplified from the general structure that minimizes the Bayes risk. And, perhaps not surprisingly, it is found that the detector/estimator resulting from the minimization of the Bayes risk consists of an infinite bank of linear correlators (matched filters) whose output, as a continuous function of the components of the hypothesized signal state vector $x$, is an essential component of the likelihood ratio. We show how the geometric properties of this likelihood ratio field over $x$ allow one to derive a simple expression for $T$ in the asymptotic limit of high signal-to-noise ratio (SNR). This asymptotic expression, valid in many applications, eliminates much of the burden of working in the continuous state parameter domain, and the simplified processing yields a concise expression for the detection criterion and a simultaneous estimate for the signal state.

In addition to simplifying the decision process, this paper produces performance predictions in the form of detection and false alarm probabilities as functions over the signal parameter domain $x$. These quantities are densities in $x$ for the probabilistic occurrence of maxima at $x$ in the decision statistic $T$. The detection probability, as it relates to the detection criterion, is treated as a function of two state vectors – the estimated signal location where the detection statistic assumes its maximum (over threshold) and the actual location of the signal. The predicted statistical discrepancy between the two provides the state estimation accuracy.

Evaluation of the false alarm density over the domain $x$ is, in general, a problem of fundamental difficulty – one that affords no reasonably usable expression that is directly obtainable. However, in our application to additive Gaussian noise and high declaration threshold, we are able to evaluate detection probabilities directly from the geometric properties we establish over the signal parameter domain. We then use stationarity relations obtained between signal and false alarm probabilities to solve the difficult problem of obtaining a concise analytical expression for the false alarm density over the domain $x$.

The expressions we obtain for the detection and false alarm probabilities over the signal parameter domain $x$ at high threshold in Gaussian noise reveal the following behavior. For a fixed decision threshold the false alarm density and detection probability are strong functions of SNR; the detection probability can range from near zero to near unity. The false alarm density is small when the detection probability is close to one or zero, and greater by orders of magnitude when the detection probability is one-half. When the SNR varies over the signal state domain $x$ by a factor well in excess of unity, regions of nearly certain and nearly impossible detection are separated by regions of marginal detection that contain most of the false alarms. The detection and false alarm descriptors we obtain address the primary concerns of a practical application – for what signal characteristics (what state) is a signal detectable and how severe is the false alarm problem.

The paper is organized as follows. Section 2 introduces the notation and derives the Bayes decision rule for joint detection and estimation. A decision statistic $T$ follows from an optimization criterion based on minimizing the Bayes risk. The decision statistic $T$ is used in Section 3 to establish an important relation between detection and false alarm probabilities by means of a general stationarity property. Section 4 then specializes to Gaussian additive noise, and the decision statistic is used again in a form that explicitly introduces the likelihood ratio as a field over the signal parameter domain $x$. Through an asymptotic approximation to the decision statistic at high threshold, the geometry of the decision statistic over the continuum of the domain $x$ is developed. This geometry is key to the analysis in Section 5 of the detector/estimator performance as a function of $x$, particularly when the SNR varies with $x$. Section 6 applies Monte Carlo simulation to a simple signal and measurement model to test the theoretical expressions of Section 5. State accuracy and false alarm occurrence are seen to agree very well with the predictions. Section 7 sums the performance over the domain $x$ to form a global operating characteristic that is used to evaluate a global false alarm probability – a quantity that would be very tedious to obtain by other means.

Significant new results are labeled as theorems.

## 2. FORMULATION

The measurement vector $\mathbf{m}$ consists of noise $\mathbf{n}$ and possibly an additive signal $\mathbf{s}(x)$ described by an unknown parameter $x$,

$$\mathbf{m} = \mathbf{s}(x) + \mathbf{n} \quad ; \quad \begin{cases} \mathbf{s}(x) = \mathbf{0} \text{ if } x = x_0 \\ \mathbf{s}(x) \neq \mathbf{0} \text{ if } x \in \chi \end{cases} . \tag{1}$$

The parameter space is the union $\{x_0\} \cup \chi$, where $x_0$ actually represents the absence of signal and the set $\chi$ is a continuous domain of dimension $d$, which in typical applications may include such dimensions as signal appearance time, position, velocity, amplitude, shape parameters, etc. The measurement or observation vector space is of finite dimension and is assumed to be continuous in its elements, admitting a probability density in $\mathbf{m}$ for all values of $x$. In a typical application the elements would include samples over space, time, and possibly frequency, etc. The problem is to decide whether signal is present, and if so, with what value of $x$. The measurement model provides the conditional probability for $\mathbf{m}$ with signal absent,

$$p_0(\mathbf{m}) \quad \text{if} \quad x = x_0 , \tag{2}$$

and the conditional probability for $\mathbf{m}$ given signal present at $x$,

$$p_1(\mathbf{m}|x) \quad \text{if} \quad x \in \chi . \tag{3}$$

These two cases are presumed to occur respectively with *a-priori* probability $a_0$ and density $a(x)$, where

$$a_0 + \int_{x \in \chi} a(x)\, dx = 1 . \tag{4}$$

It should be noted that the general properties developed in this section and the next are good for a more general measurement model, in which the noise randomly alters the signal; the only assumption is that a useful measurement yields values of the likelihood ratio

$$L(\mathbf{m}|x) = p_1(\mathbf{m}|x) / p_0(\mathbf{m}) \tag{5}$$

significantly exceeding unity for probable values of $\mathbf{m}$ when signal is present. The additive form (1) is intended for the analysis of Gaussian noise in later sections.

The unconditional probability density for $\mathbf{m}$ is

$$p(\mathbf{m}) = a_0 p_0(\mathbf{m}) + \int_{x \in \chi} a(x) p_1(\mathbf{m}|x)\, dx \tag{6}$$

It follows that the Bayes probabilities for no signal $(x = x_0)$ or for signal at $x \in \chi$, conditioned on a particular measurement $\mathbf{m}$, are

$$b_0(\mathbf{m}) = a_0 p_0(\mathbf{m}) / p(\mathbf{m}),$$
$$b(x|\mathbf{m}) = a(x) p_1(\mathbf{m}|x) / p(\mathbf{m}); \tag{7}$$

that is, the random variable $x_0$ occurs with probability $b_0(\mathbf{m})$, and the random variable $x \in \chi$ occurs with conditional probability density $b(x|\mathbf{m})$. These are the constituents of an optimum detection/estimation rule, to be defined and assembled from a cost matrix as follows [8].

The standard Bayesian approach to assigning costs to decisions is appropriate with the notation extended to accommodate different costs at different signal locations $x$ and declaration locations $x'$. We let the first argument of the cost function $C(x', x)$ represent the declaration location and the second the actual signal location. Furthermore, we represent the signal-absent parameter $x_0$ explicitly and assume the location parameters $x, x'$ are otherwise from the set $\chi$. Accordingly we have:

$C(x', x)$:   Estimation error: cost for a detection at $x'$ of a signal at $x$

$C(x_0, x)$:   Missed detection: cost for declaring signal absent when it is present with parameter $x$

$C(x', x_0)$:   False alarm: cost for a signal declaration at $x'$ when signal is absent

$C(x_0, x_0)$:   Cost for declaring signal absent when it is absent  .

The cost function is assumed to be continuous and bounded in $x$ and $x'$, and $C(x', x)$ is assumed to attain a minimum over $x'$ at $x' = x$, the correct state estimate. The notation becomes simplified by setting $C(x_0, x_0)$ to zero (without loss of generality), by redefining the false-alarm cost $C(x', x_0)$ as $C(x')$, and by defining a "reward" function $R(x', x) = C(x_0, x) - C(x', x)$ as the cost difference between declaring a signal at $x$ to be absent and declaring it to be present at $x'$. We assume a reasonable cost assignment that would leave $R(x', x)$ nonnegative but at or near zero for $x'$ far from $x$. The Bayes risk is defined as the cost function averaged over the Bayesian priors and the measurement distribution:

$$B_{risk}(x', M_C) = \int_{\{\mathbf{m} \in M_C, x'\}} \int C(x', x) a(x) p_1(\mathbf{m}|x) dx d\mathbf{m}$$
$$+ \int_{\{\mathbf{m} \notin M_C, x'\}} \int C(x_0, x) a(x) p_1(\mathbf{m}|x) dx d\mathbf{m} \tag{8}$$
$$+ \int_{\{\mathbf{m} \in M_C, x'\}} C(x') a_0 p_0(\mathbf{m}) d\mathbf{m} \quad.$$

The integration domain $\{\mathbf{m} \in M_C, x'\}$ represents the critical region of measurement space that is assigned a signal-present decision at $x'$. Measurement values $\mathbf{m}$ that lie in the complementary region $\{\mathbf{m} \notin M_C, x'\}$ receive a no-signal decision.

Our optimization criterion for simultaneous detection and estimation can now be stated explicitly: *selection of $x'$ and the corresponding critical region $\{\mathbf{m} \in M_C\}$ to minimize the Bayes risk*. This optimization balances the cost of a missed detection, the cost of misestimating the signal state for a valid detection, and the cost of a false detection. The decision structure that satisfies this criterion is readily obtained, and in anticipation we define the decision statistic

$$T(x'|\mathbf{m}) = \frac{1}{p(\mathbf{m})} \left[ \int R(x', x) p_1(\mathbf{m}|x) a(x) dx - C(x') p_0(\mathbf{m}) a_0 \right] \tag{9}$$

The decision statistic $T(x'|\mathbf{m})$ is readily seen to be the difference between the posterior loss for a signal-absent decision and that for a signal-at-$x'$ decision conditional on $\mathbf{m}$:

$$T(x'|\mathbf{m}) = \int \left[ C(x_0, x) - C(x', x) \right] b(x|\mathbf{m}) dx$$
$$+ \left[ C(x_0, x_0) - C(x', x_0) \right] b_0(\mathbf{m}) \quad. \tag{10}$$

With the recognition that the measurement probabilities are normalized to integrate to unity over the combined region, $\{\mathbf{m} \in M_C, x'\} \cup \{\mathbf{m} \notin M_C, x'\}$, the Bayesian risk can be represented as

$$B_{risk}(x', M_C) = \int C(x_0, x) a(x) dx$$
$$- \int_{\{\mathbf{m} \in M_C, x'\}} T(x'|\mathbf{m}) p(\mathbf{m}) d\mathbf{m}. \tag{11}$$

The equation above shows clearly that choosing $x'$ and $M_C$ to maximize the integration of $T(x'|\mathbf{m})$ over $\mathbf{m}$ minimizes the Bayes risk, and it is very obvious how the choice should be made: values of $\mathbf{m}$ that yield a negative $T(x'|\mathbf{m})$ for all $x'$ should be excluded from the critical region $M_C$, and $x'$ should be chosen to maximize the positive values that are included in $M_C$. The definition of $T(x'|\mathbf{m})$ ensures that it is negative when a signal-absent decision minimizes the expected cost for the measurement $\mathbf{m}$ and positive for some $x'$ when a signal-present decision minimizes the cost. Thus the combined detection-estimation decision finds the signal parameter estimate $x' = x_\mathbf{m}$ that maximizes $T(x'|\mathbf{m})$ (minimizes the net expected cost),

$$T(x_{\mathbf{m}}|\mathbf{m}) = \max_{x'}\{T(x'|\mathbf{m})\} \quad , \qquad (12)$$

and then, as the decision test, chooses between signal present and signal absent:

$$T\left(x_{\mathbf{m}}|\mathbf{m}\right) \quad \begin{cases} > 0\,, & \text{choose } x_{\mathbf{m}} \\ < 0\,, & \text{choose } x_0 \end{cases} \quad , \qquad (13)$$

resulting in the critical region

$$\{\mathbf{m} \in M_C, x'\} = \{\mathbf{m} \in T\left(x_{\mathbf{m}}|\mathbf{m}\right) > 0, x' = x_{\mathbf{m}}\} \quad . \qquad (14)$$

The development to this point has produced a statistic that provides the proper choice between signal and no-signal and simultaneously specifies the signal parameter estimate $x_{\mathbf{m}}$ that best explains a signal decision. However, analysis of the decision statistic can provide much more information. It has the potential to describe the distribution of signal location about the observed maximum at $x' = x_{\mathbf{m}}$, and it can be used in a predictive sense (with $\mathbf{m}$ as a random variable) to determine signal detection and false alarm probabilities over the signal parameter domain. In the next section we manipulate the form of the decision statistic to produce general expressions for these probabilities. In the succeeding sections on Gaussian noise, we develop explicit analytic expressions for these probabilities (at high SNR, high threshold).

## 3. GENERAL PROPERTIES

Our objective is to mold the decision statistic into a form that reveals performance over the domain $x$ of signal parameters. The following basic operations are required. We adopt as our general measure of performance,

$$\overline{T} = \int_{T\left(x_{\mathbf{m}}|\mathbf{m}\right)>0} T\left(x_{\mathbf{m}}|\mathbf{m}\right) p(\mathbf{m}) d\mathbf{m} \quad . \qquad (15)$$

It is obvious from (11) that $\overline{T}$ represents the maximum obtainable reduction in the Bayes risk from $\int C\left(x_0, x\right) a\left(x\right) d\,$, the expected cost of a missed detection as calculated from the Bayes prior alone. The declaration threshold, which depends on reward and cost functions, appears as a limit of integration, and to examine dependence carried by this threshold, we introduce a scale factor into the cost function. The performance measure $\overline{T}$, now dependent on a scaled false-alarm cost, is converted to the integral of a density over the parameter space $x'$, and the expression reveals quantities that can be identified as probability densities for signal detection and false alarm occurrence. Finally the sensitivity of the performance measure density to the scaled cost function provides a relation between the detection and false alarm probabilities over all points $x'$ of the parameter domain. While the expressions for these probabilities are too complex for direct evaluation, we shall find we can determine explicit analytic forms for them under the assumptions and approximations of the Gaussian analysis.

We begin the development sketched above by replacing $C(x')$ with a scaled cost function $\mu C(x')$, where the scalar factor $\mu$ will provide a mechanism for the construction of an operating characteristic. We denote the decision statistic of (9) corresponding to the scale factor $\mu$ as $T(x', \mu|\mathbf{m})$. Consider now, as the performance measure, the decision statistic weighted with respect to the marginal probability density $p(\mathbf{m})$ over all measurement realizations $\mathbf{m}$ that produce a signal decision ($T(x_{\mathbf{m}}, \mu|\mathbf{m}) > 0$):

$$\overline{T}\left(\mu\right) = \int_{T(x_{\mathbf{m}}, \mu|\mathbf{m})>0} T(x_{\mathbf{m}}, \mu|\mathbf{m}) p(\mathbf{m}) d\mathbf{m} \quad . \qquad (16)$$

Sensitivity of the performance model to the cost factor $\mu$ can be immediately obtained from the expression in (9) (with $C(x')$ replaced by $\mu C(x')$) as

$$\frac{\partial}{\partial \mu}\overline{T}\left(\mu\right) = -\int_{T(x_{\mathbf{m}}, \mu|\mathbf{m})>0} C\left(x_{\mathbf{m}}\right) p_0\left(\mathbf{m}\right) a_0 d\mathbf{m} \quad . \qquad (17)$$

The integration region depends on $\mu$, but its variation does not contribute to the expression, because the integrand is zero at the threshold $T\left(x_{\mathbf{m}}, \mu|\mathbf{m}\right) = 0$.

We wish to distribute the averaged decision statistic over the parameter space of potential declaration locations $x'$ as the first step toward introducing signal-detection and false-alarm densities over that space. The Dirac delta function $\delta\left(x'\right)$, with the property that

$$\int \delta\left(x'\right) f\left(x'\right) dx' = f\left(0\right) \qquad (18)$$

for any suitably continuous function $f\left(x'\right)$ [9], serves as a convenient vehicle for this partitioning. The delta function can be inserted into (16) as

$$\overline{T}\left(\mu\right) = \int_{T(x_{\mathbf{m}}, \mu|\mathbf{m})>0} \left\{ \int \delta\left(x' - x_{\mathbf{m}}\right) T(x', \mu|\mathbf{m}) dx' \right\} p(\mathbf{m}) d\mathbf{m} \quad , \qquad (19)$$

and an interchange of the order of integration yields the desired density form:

$$\overline{T}\left(\mu\right) = \int \tau\left(x'|\mu\right) dx' \qquad (20)$$

where

$$\tau\left(x'|\mu\right) = \int_{T(x_{\mathbf{m}}, \mu|\mathbf{m})>0} \delta\left(x' - x_{\mathbf{m}}\right) T(x', \mu|\mathbf{m}) p(\mathbf{m}) d\mathbf{m} \quad . \qquad (21)$$

Substitution of the expression for $T(x',\mu|\mathbf{m})$ from (9) enables the introduction of the signal-detection and false-alarm densities we seek:

$$\tau(x'|\mu) = \int R(x',x)\,p_D(x'|x,\mu)\,a(x)\,dx \\ -\mu C(x')\nu_F(x'|\mu)\,a_0 \quad, \tag{22}$$

where

$$p_D(x'|x,\mu) = \int_{T(x_\mathbf{m},\mu|\mathbf{m})>0}\delta(x'-x_\mathbf{m})\,p_1(\mathbf{m}|x)\,d\mathbf{m} \\ \nu_F(x'|\mu) = \int_{T(x_\mathbf{m},\mu|\mathbf{m})>0}\delta(x'-x_\mathbf{m})\,p_0(\mathbf{m})\,d\mathbf{m} \quad. \tag{23}$$

Here $p_D(x'|x,\mu)$ can be recognized as the probability density, in $x'$, for detection at $x'$ of a signal which actually occurred at $x$, and $\nu_F(x'|\mu)$ as the density of false detections at $x'$. The delta functions provide the dimensionality for densities in $x'$, and the integration over $\mathbf{m}$ sums the probabilities for which $x_\mathbf{m} = x'$.

The sensitivity of the density $\tau(x'|\mu)$ to the cost factor $\mu$ can now be used to generate the relation between signal detection density and false-alarm density that we seek. Differentiation of the expression in (21) yields

$$\frac{\partial}{\partial\mu}\tau(x'|\mu) = \int_{T(x_\mathbf{m},\mu|\mathbf{m})>0}\delta(x'-x_\mathbf{m})\frac{\partial}{\partial\mu}T(x',\mu|\mathbf{m})\,p(\mathbf{m})\,d\mathbf{m} \\ = -\int_{T(x_\mathbf{m},\mu|\mathbf{m})>0}\delta(x'-x_\mathbf{m})\,C(x')\,p_0(\mathbf{m})\,a_0\,d\mathbf{m} \\ = -C(x')\nu_F(x'|\mu)\,a_0 \quad. \tag{24}$$

As in the derivation of (17), the variation in the threshold does not contribute because the integrand is zero there. Differentiation of the same quantity from the expression in (22) yields

$$\frac{\partial}{\partial\mu}\tau(x'|\mu) = \int R(x',x)\frac{\partial}{\partial\mu}p_D(x'|x,\mu)\,a(x)\,dx \\ -\mu C(x')\frac{\partial}{\partial\mu}\nu_F(x'|\mu)\,a_0 - C(x')\nu_F(x'|\mu)\,a_0 \quad. \tag{25}$$

Comparison with (24) shows that the first two terms on the right sum to zero. The combination of these two terms can be identified as the variation of $\tau(x'|\mu)$ in (21) with respect to the change in integration domain, and this change occurs at the threshold where the integrand itself is zero and results in the cancellation of the terms. Now the variational relation we seek follows from this equality of the first two terms and takes the form of the following differential equation:

$$\int R(x',x)\frac{\partial}{\partial\mu}p_D(x'|x,\mu)\,a(x)\,dx = \mu C(x')\frac{\partial}{\partial\mu}\nu_F(x'|\mu)\,a_0 \quad. \tag{26}$$

The detection and false-alarm densities $p_D(x'|x,\mu)$ and $\nu_F(x'|\mu)$ are densities in the signal parameter space. Differentiation with respect to $\mu$ renders the derivatives as densities with respect to amplitude as well. We shall see when we analyze the high-SNR, Gaussian-noise case in Section 5 that the detection probability can be approximated in a reasonably direct manner. However, the false alarm probability involves the probabilistic occurrence of local maxima over the field $x'$ of possible signal parameters, and direct evaluation is not so readily attained. For false alarms we shall use the differential relation (26) to evaluate the false alarm probability from the detection probability, as shown in general in Theorem 1.

*Theorem 1*: The stationarity relation of (26) can be integrated from an upper limit $\mu = \infty$, where the infinite false alarm cost eliminates all target declarations and where both detections and false alarms vanish, to $\mu = 1$, where the assigned cost function $C(x')$ applies, to give

$$\nu_F(x'|1) = -\int_1^\infty\frac{1}{\mu}\frac{\partial}{\partial\mu}\int\frac{R(x',x)}{C(x')}p_D(x'|x,\mu)\frac{a(x)}{a_0}\,dx\,d\mu. \tag{27}$$

If the lower integration limit is set at $\mu$ instead of unity, the relation above serves as a generalized version of the local operating characteristic of the detector as a function of arbitrary threshold (provided $\mu C(x')$ is taken as the cost function determining the threshold).

These relations bear a resemblance to those for the operating characteristic for the binary decision problem addressed in, for example, Van Trees [10] and Schreier and Scharf [11], with $\mu$ corresponding to a likelihood threshold. However, in the binary problem the operating characteristic is obtained very directly from the likelihood ratio threshold that results from both the Neyman-Pearson criterion and a two-possibility Bayes criterion. In our analysis, on the other hand, the relation between false alarm and detection probabilities is obtained only indirectly from the basic detection criterion of Equations (9), (12), and (13). Our object of primary interest is the detection and false alarm behavior over the signal parameter space, and the basic detection criterion only implicitly comprehends the dependence on $x'$. The conversion of the detection criterion to an expression that includes probability densities over $x'$ incorporates the threshold dependence required to balance rewards and costs over the entire space. Scaling the false-alarm cost function by $\mu$ provides the variation in the performance measure that yields a stationarity relation connecting signal-detection and false-alarm probabilities.

## 4. ADDITIVE GAUSSIAN NOISE

The specific case of Gaussian noise additive on a family of known signals is both widely useful and readily analyzable. The immediate aim of the following analysis is to construct an approximation, good in the asymptotic limit of high SNR, to the decision statistic $T(x'|\mathbf{m})$ in (9). In the course of this analysis a number of interesting properties of $T(x'|\mathbf{m})$ will emerge, including the intrinsic estimation accuracy and estimates of both local and global detection performance over $x'$ in the form of detection and false alarm probabilities.

Under the conditions we impose for the analysis, the parameter estimate which maximizes the decision statistic is essentially the same as that which maximizes the likelihood ratio over the same local region. In fact, the likelihood ratio of (5) appears explicitly in the decision statistic if we move $p_0(\mathbf{m})$ outside the bracket in (9):

$$T(x'|\mathbf{m}) = \frac{p_0(\mathbf{m})}{p(\mathbf{m})} \left[ \int R(x',x) L(\mathbf{m}|x) a(x)\, dx - C(x') a_0 \right]. \quad (28)$$

We begin by exploring the properties of the log-likelihood ratio (LLR). The measurement model of Equation (1) applies with a noise vector $\mathbf{n}$ of mean zero and covariance $\mathbf{C}$. In the Gaussian case the likelihood ratio of (5) has as its denominator, the Gaussian noise probability of $\mathbf{n}$, and as its numerator, the Gaussian noise probability of $\mathbf{m}$-$\mathbf{s}$. The resulting expression for the LLR is

$$\lambda(\mathbf{m}|x) = \log L(\mathbf{m}|x)$$
$$= \mathbf{s}*(x)\mathbf{C}^{-1}\mathbf{m} - \frac{1}{2}\mathbf{s}*(x)\mathbf{C}^{-1}\mathbf{s}(x), \quad (29)$$

where the asterisk denotes vector transpose. We note that the first term, representing the output of a linear matched filter for the hypothesized signal applied to the measurement, has a Gaussian probability distribution with mean zero in the absence of signal and with mean $\mathbf{s}*(x)\mathbf{C}^{-1}\mathbf{s}(x)$ when the particular hypothesized signal is present. The variance of this term is $r^2(x) = \mathbf{s}*(x)\mathbf{C}^{-1}\mathbf{s}(x)$, where $r(x)$ can be recognized as the amplitude SNR for $\mathbf{s}(x)$. Defining the matched-filter output normalized to unit variance (with signal either present or absent),

$$y(\mathbf{m}|x) = \mathbf{s}*(x)\mathbf{C}^{-1}\mathbf{m} / r(x), \quad (30)$$

we can write

$$\lambda(\mathbf{m}|x) = r(x) y(\mathbf{m}|x) - \frac{1}{2} r^2(x) \quad (31)$$

and, when the hypothesized signal matches the actual signal,

$$\langle y(\mathbf{m}|x) \rangle = \begin{cases} r(x), & \text{signal present at } x \\ 0, & \text{signal absent} \end{cases}. \quad (32)$$

The expectation operator $\langle \cdot \rangle$, here and subsequently, represents the statistical expectation with respect to the distribution of the noise $\mathbf{n}$. In some elementary treatments the quantity $y$ is itself used as the decision statistic [12]. This simplification is permissible for one signal but not for several competing signal hypotheses of differing SNR, where the entire form above is essential.

The quantity $y(\mathbf{m}|x)$ nevertheless plays an important intermediate part in the analysis of the LLR. The response of the wrong matched filter $\mathbf{s}(x')$ to a signal at $x$ can be computed in terms of the scaled inner product

$$\sigma(x',x) = \mathbf{s}*(x')\mathbf{C}^{-1}\mathbf{s}(x) / r(x') r(x), \quad (33)$$

with

$$\sigma(x,x) = 1,$$
$$|\sigma(x',x)| < 1, \quad x' \neq x \quad (34)$$

by the Schwarz inequality. It follows from this that

$$y(\mathbf{s}(x)|x') = r(x) \sigma(x',x) \quad (35)$$

so that the wrong filter preserves the amplitude of the signal with a decorrelation loss from the signal dissimilarity.

For notational convenience we now drop the explicit measurement dependence of $y$ and represent $y(\cdot|x')$ as $y(x')$ (and will later do the same for $\lambda(\cdot|x')$). In the absence of an actual signal in the measurement, the gradient of $y(x')$, represented here as derivatives with respect to the $d$ parameters of the hypothesized signal, $\partial y(x')/\partial x_i'$, $i = 1,2,...,d$, consists of Gaussian-distributed variables with the property that they are statistically independent of the amplitude at the hypothesized signal parameter $x'$, as seen from

$$0 = \frac{\partial}{\partial x_i'} \langle y^2(x')/2 \rangle|_N = \langle y(x') \frac{\partial}{\partial x_i'} y(x') \rangle|_N, \quad (36)$$

where the condition $N$ (noise) indicates no signal present. The components of the gradient are in general mutually correlated,

$$\langle \frac{\partial y}{\partial x_i'} \frac{\partial y}{\partial x_j'} \rangle|_N = M_{ij} = (M)_{ij}, \quad i,j = 1,2,...,d. \quad (37)$$

Similarly, the gradient of (36) gives

$$0 = M_{ij} + \left\langle y(x')\frac{\partial^2}{\partial x_i'\partial x_j'}y(x')\right\rangle\Big|_N , \qquad (38)$$

so that the matrix of second derivatives is negatively correlated with amplitude. More particularly, if we introduce the matrix $U$ as the difference between the second gradient of $y(x')$ and what happens to be its conditional mean (conditioned by amplitude $y(x')$),

$$\frac{\partial^2 y}{\partial x_i'\partial x_j'} = -yM_{ij} + U_{ij}, \qquad (39)$$

we can use the basic characteristics we infer about $U$ to provide some general geometric properties of the field of normalized matched-filter outputs. Since $y(x')$ and all its derivatives are mean zero in the absence of signal, taking the expectation of (39) as the equation stands and after multiplication by $y(x')$ (and using (38) and the normalization condition $\langle y^2(x')\rangle|_N = 1$) yields

$$\langle U\rangle|_N = 0, \ \ \langle yU\rangle|_N = 0. \qquad (40)$$

This signifies that, for an observed value $y(x')$, the second gradient has a mean part proportional to $y(x')$ times what we may call the *curvature matrix M,* plus a random Gaussian part $U$ independent of amplitude. As $y(x')$ increases, the random part becomes asymptotically unimportant and the curvature becomes more accurately described by $M$. Thus as we survey the field of normalized matched-filter outputs over the signal parameter domain, an ambient peak (signal absent) that might trigger a detection declaration in the limit of high threshold (high false alarm cost $C(x')$) would have its shape characterized by the curvature matrix $M$. Note also that (30) and (38) imply

$$\mathrm{cov}\big(y(x'), y(x)\big) = \sigma(x', x) \qquad (41)$$

and

$$\begin{aligned}\frac{\partial^2}{\partial x_i'\partial x_j'}\sigma(x',x)\Big|_{x'=x} &= \left\langle\left[\frac{\partial^2}{\partial x_i'\partial x_j'}y(x')\Big|_{x'=x}\right]y(x)\right\rangle\Big|_N \\ &= -M_{ij}(x) .\end{aligned} \qquad (42)$$

As one might expect, the shape of a local peak in the matched-filter output provides no additional information about the signal's presence; from (35) a signal at $x$ with amplitude $r(x)$ produces a mean curvature

$$\left\langle\frac{\partial^2}{\partial x_i'\partial x_j'}y(x')\Big|_{x'=x}\right\rangle\Big|_{s(x)} = -r(x)M_{ij}(x) \qquad (43)$$

indistinguishable from a spontaneous noise peak having a comparable amplitude $y \approx r$ .

The foregoing enables us to describe the statistical geometry of maxima in the LLR. Reinserting the scaling in (31), we find that

$$\frac{\partial}{\partial x_i'}\lambda(x') = r(x')\frac{\partial}{\partial x_i'}y(x') + \big(y(x')-r(x')\big)\frac{\partial}{\partial x_i'}r(x') \qquad (44)$$

and at a maximum, $\partial\lambda/\partial x_i' = 0$ ,

$$\frac{\partial y}{\partial x_i'} = \left(1-\frac{y}{r}\right)\frac{\partial r}{\partial x_i'}, \qquad (45)$$

so that variable SNR implies that maxima in $y(x')$ and $\lambda(x')$ do not quite coincide. With this constraint on $\partial y/\partial x_i'$, the matrix of second gradients of $\lambda$ becomes

$$\begin{aligned}\frac{\partial^2\lambda}{\partial x_i'\partial x_j'} &= r\frac{\partial^2 y}{\partial x_i'\partial x_j'} - \frac{\partial r}{\partial x_i'}\frac{\partial r}{\partial x_j'} \\ &+ \left(1-\frac{y}{r}\right)\left(2\frac{\partial r}{\partial x_i'}\frac{\partial r}{\partial x_j'} - r\frac{\partial^2 r}{\partial x_i'\partial x_j'}\right) .\end{aligned} \qquad (46)$$

We now argue that the last term of (46) can be neglected because $1-y/r$ will be small for a declarable peak in a high SNR region of the signal parameter domain, i.e., that $y(x')$ will be large and near $r(x')$ in magnitude at a local maximum in log-likelihood – whether or not a signal is present. A detailed explanation follows.

For signal present it has been established in (32) that $\langle y\rangle = r$ for $y$ evaluated at a hypothesized location $x'$ that matches the location $x$ of an actual signal. Our analysis is restricted to high thresholds where only signals residing in regions of high SNR, $r^2 >> 1$, can be expected to be detected. The quantity $y$ is a normalized Gaussian variable measured in units of standard deviations. A difference of many standard deviations at the signal location is very unlikely, and so, with high probability, $1-y/r$ will be small. When no signal is present ($\langle y\rangle = 0$), it is also very unlikely that $\lambda = y^2/2 - (y-r)^2/2$ (Equation (31)) will be large at any particular location $x'$. However, a large number of effectively "independent regions" will be searched over the domain of signal parameters, and the probability of encountering a false alarm with a high log-likelihood value may be substantial. If such a false alarm occurs, $y$ must be large, since $\lambda \le y^2 / 2$. The second term of $\lambda$, $(y - r)^2 / 2$, is subtracted, and if it is not small (if $y$ is not near $r$ in value), achieving a given level of $\lambda$ would require an even larger and less probable value of $y$. Thus in a parameter domain in which $r$ varies significantly over $x$, the most probable configuration for achieving a high specific (threshold) log-likelihood value is with $y \approx \sqrt{2\lambda}$ and therefore $r(x') \approx y$.

This brief analysis, by the way, provides an intuitive glimpse of one of the findings of the paper that is derived more formally later: that false alarms tend to concentrate in regions of the domain $x$ of marginal detectability.

With the last term of (46) omitted and $y$ replaced by $r$ as the coefficient of the curvature matrix (39), the second gradient of LLR at a local maximum in the signal parameter field becomes

$$\frac{\partial^2}{\partial x_i' \partial x_j'} \lambda(x') \cong -r^2(x') Q_{ij}(x') ,\qquad (47)$$

where

$$Q_{ij} = M_{ij} + (\frac{\partial}{\partial x_i'} \log r)(\frac{\partial}{\partial x_j'} \log r) .\qquad (48)$$

Note that this equation is exact for the expected value of $\partial^2 \lambda / \partial x_i' \partial x_j'$ when the signal is present at $x' = x$. The matrix describing the LLR maximum shows an increased curvature along the axis of the local SNR gradient; the measured SNR provides additional localization information in that direction.

Now we return to evaluating the decision statistic. Approximating the integral in $T(x'|\mathbf{m})$ (28) requires us to specify the function $R(x',x)$, and particularly the tolerance for estimation error. To do this meaningfully we note that a likelihood peak has intrinsic dimensions, as seen in the asymptotic approximation obtained by retaining the first two terms in the Taylor expansion of $\lambda$ around the location $x_{\mathbf{m}}$ of the maximum $\lambda_{\mathbf{m}} = \lambda(\mathbf{m}|x_{\mathbf{m}})$ (which, under the conditions of our analysis, is essentially the same as the location of the corresponding local maximum in $T(x'|\mathbf{m})$):

$$L(\mathbf{m}|x) \cong \exp \left[ \lambda_{\mathbf{m}} - \frac{r^2(x_{\mathbf{m}})}{2}(x - x_{\mathbf{m}}) * Q(x_{\mathbf{m}})(x - x_{\mathbf{m}}) \right].\qquad (49)$$

This approximation, good for large $r^2$, indicates that most of the likelihood occurs within a region $X(x_{\mathbf{m}})$ defined by

$$\frac{1}{2}(x - x_{\mathbf{m}}) * Q(x_{\mathbf{m}})(x - x_{\mathbf{m}}) \le c / r^2 \qquad (50)$$

with $c$ being a modest constant; e.g., for 95 percent likelihood containment, $c = 3.0, 3.8, 4.7$ for dimension $d = 2, 3, 4$, respectively. One other assumption we shall adopt throughout is that the signal prior $a(x)$, the reward function $R(x', x)$, the cost function $C(x')$, and the curvatures are slowly varying on the scale of the peak dimensions.

This characterization shows that, for our high-SNR Gaussian application, the reward function $R(x', x)$, which in general enforces the requirement that the estimated location be near the true signal location, has little quantitative impact provided that it extends at significant amplitude over the region of a peak. In that case it becomes convenient to define its region of significance to match the extent of the peak:

$$R(x', x) = \begin{cases} 1, & x \in X(x') \\ 0, & \text{otherwise} \end{cases}.\qquad (51)$$

In choosing unity above instead of a variable $R(x', x')$, we recognize that no generality is lost, because the results to follow depend only on the ratio $R(x', x') / C(x')$. Then for integration over the domain $x \in X(x_{\mathbf{m}})$ we have, for the integral in (28), evaluated at $x' = x_{\mathbf{m}}$,

$$\int R(x_{\mathbf{m}}, x) a(x) L(\mathbf{m}|x) dx = e^{\lambda_{\mathbf{m}}} r^{-d}(x_{\mathbf{m}}) a(x_{\mathbf{m}}) V(x_{\mathbf{m}}) \qquad (52)$$

where

$$V(x_{\mathbf{m}}) = \left\| Q(x_{\mathbf{m}}) / 2\pi \right\|^{-1/2} \qquad (53)$$

defines a $d$-dimensional volume scale associated with the local likelihood-ratio peak. It provides a measure of the local region in signal state space over which the matched-filter output would be affected by the presence of a signal, and it also represents the region across which significant noise correlations would exist. Thus a large volume scale indicates few independent noise realizations (separate false alarm possibilities) in the vicinity of $x_{\mathbf{m}}$ and correspondingly boosts the decision statistic through the increase in the term (52); a small volume scale adjusts the decision statistic downwards as compensation for a greater density of effectively independent noise draws.

## 5. PERFORMANCE

We are now equipped to evaluate the form the decision criterion assumes as a threshold on LLR, to calculate the estimation error in signal parameter space for a detection declaration, and to calculate detection and false-alarm probabilities, all in the high-SNR limit of our Gaussian application.

The decision test at high thresholds in additive Gaussian noise can be greatly simplified from the full processing indicated by Equations (9), (12), and (13). The direct use of the Bayesian decision rule involves a multidimensional integration over the domain $x$ for every hypothesized signal location $x'$ – a formidable undertaking – followed by seeking the maximum. But in the Gaussian case in the high-SNR limit, the integrations (on the decision statistic in the form (28)) can be performed analytically, as shown in (52). The search for $x_{\mathbf{m}}$ must still be carried out, but it is accomplished by forming the LLR over the domain of $x$, finding the local maxima, and comparing them to the local thresholds. If it happens that more than one local maximum is over threshold (and if that matters in an application), it is a simple matter to check the factors of the threshold to determine which local maximum corresponds to the global maximum $x_{\mathbf{m}}$ of the original Bayes rule. The decision rule can be converted to a test of the LLR against a threshold – a threshold that varies over the signal parameter domain through the dependence of prior probabilities, the cost function, and the fundamental geometry of the likelihood peaks. The use of (52) in (28) produces the quantitative prescription of Theorem 2.

*Theorem 2:* The application of the decision test (13) in the high-SNR limit of our Gaussian application $\left( r^2 >> 1 \right)$ takes the form,

$$\lambda_{\mathbf{m}} \begin{cases} > \lambda_r(x_{\mathbf{m}}) \Rightarrow \text{choose signal at } x_{\mathbf{m}}, \\ < \lambda_r(x_{\mathbf{m}}) \Rightarrow \text{choose no-signal} \end{cases},$$
$$\lambda_r(x') = \log\left[ a_0 C(x') r^d(x') / a(x') V(x') \right]. \tag{54}$$

The threshold $\lambda_r(x')$ does not depend on the measurement $\mathbf{m}$ and can be computed anywhere over the signal parameter domain, but it need be evaluated only where the LLR has a local maximum: $x' = x_{\mathbf{m}}$. This equation constitutes a major result of this paper, furnishing an asymptotically correct, simple implementation of the Bayes solution for detection and estimation. The variable threshold $\lambda_r(x')$ contains the entire prior information and geometric characteristics ($V(x')$) of the signal parameter domain as well as a modest dependence on SNR. Note particularly that the scale volume $V(x')$, which depends on the choice of the coordinate system $x$, occurs in the combination $a(x')V(x')$, which is invariant to the coordinate choice and represents the prior signal probability encompassed by the volume.

The accuracy of the state estimate $x_{\mathbf{m}}$ can be readily analyzed through a Taylor expansion about the true signal location $x_{\mathbf{s}}$. Use of quantities defined in Section 4 in the expansion provides a statistical representation of both the signal state estimate and the peak LLR amplitude. Equations (33) to (35) show directly that $\left\langle y(x') \right\rangle\big|_{s(x_s)}$ has a local maximum at $x' = x_{\mathbf{s}}$. It then follows immediately by taking expectations of the equation for the LLR gradient (44) that $\left\langle \lambda(x') \right\rangle\big|_{s(x_s)}$ also has a zero gradient and hence a local maximum at $x' = x_{\mathbf{s}}$. Although the expected value of the LLR has a local maximum at $x' = x_{\mathbf{s}}$, the peak LLR will be displaced by nonzero noise-induced slope (44), whose covariance is [see (36) and (37)]

$$\left\langle \frac{\partial \lambda}{\partial x_i'}(x_{\mathbf{s}}) \frac{\partial \lambda}{\partial x_j'}(x_{\mathbf{s}}) \right\rangle\Big|_{s(x_s)} = r^2(x_{\mathbf{s}}) M_{ij}(x_{\mathbf{s}}) + \frac{\partial r}{\partial x_i'}(x_{\mathbf{s}}) \frac{\partial r}{\partial x_j'}(x_{\mathbf{s}})$$
$$= r^2(x_{\mathbf{s}}) Q_{ij}(x_{\mathbf{s}}). \tag{55}$$

Denoting the LLR gradient vector at the signal location as

$$\nabla \lambda_{\mathbf{s}} = \left( \frac{\partial \lambda}{\partial x_1'}(x_{\mathbf{s}}), \dots, \frac{\partial \lambda}{\partial x_d'}(x_{\mathbf{s}}) \right)^T \tag{56}$$

and asymptotically replacing the curvature by its average in the high SNR limit [see (47)], the expansion about the signal location (where $\lambda(x_{\mathbf{s}})$ is denoted $\lambda_{\mathbf{s}}$) becomes

$$\lambda(x') = \lambda_{\mathbf{s}} + (x' - x_{\mathbf{s}}) * \nabla \lambda_{\mathbf{s}}$$
$$- \frac{r^2(x_{\mathbf{s}})}{2} (x' - x_{\mathbf{s}}) * Q(x_{\mathbf{s}})(x' - x_{\mathbf{s}}). \tag{57}$$

Taking the gradient of (57) and setting it equal to zero, we see that the LLR is a maximum at the location $x_{\mathbf{m}}$ (unbiased under our approximations and assumptions),

$$x_{\mathbf{m}} - x_{\mathbf{s}} = r^{-2} Q^{-1} \nabla \lambda_{\mathbf{s}}, \tag{58}$$

with the value

$$\lambda_{\mathbf{m}} = \lambda_{\mathbf{s}} + \frac{1}{2r^2} \nabla \lambda_{\mathbf{s}} * Q^{-1} \nabla \lambda_{\mathbf{s}}$$
$$= \lambda_{\mathbf{s}} + \frac{1}{2} (x_{\mathbf{m}} - x_{\mathbf{s}}) * r^2 Q(x_{\mathbf{m}} - x_{\mathbf{s}}). \tag{59}$$

Averaging the estimation error of (58) and the LLR peak amplitude of (59) leads to the following fundamental result.

*Theorem 3:* The estimation error in the high-SNR limit approaches the Cramer-Rao bound [13]:

$$\langle (x_{\mathbf{m}} - x_{\mathbf{s}})(x_{\mathbf{m}} - x_{\mathbf{s}})* \rangle = r^{-2}(x_{\mathbf{s}})Q^{-1}(x_{\mathbf{s}}). \tag{60}$$

Note that the *rms* estimation error in each dimension is inverse to the SNR. The expected amplitude of the LLR peak is biased upward from the expected amplitude at the signal by the second term of (59), which involves a quadratic form with the mean-zero random vector $\nabla\lambda_{\mathbf{s}}$ and its inverse covariance $Q^{-1}/r^2$. This combination produces a chi-square variable with mean $d$ and variance $2d$ [14] and yields the expected LLR peak amplitude,

$$\langle \lambda_{\mathbf{m}} \rangle = \langle \lambda_{\mathbf{s}} \rangle + d/2 = r^2(x_{\mathbf{s}})/2 + d/2. \tag{61}$$

At this point we are equipped to predict performance of the Bayes test through probabilities of detection and false alarm occurrence. The detection probability is readily obtained under our Gaussian, high-SNR assumptions through the expansion about the signal location $x_{\mathbf{s}}$ presented above. The false alarm probability, on the other hand, depends fundamentally on the occurrence of local maxima over a random field and is not so readily obtained. Equation (23) provides a theoretical representation of the probability but does not provide a feasible method for evaluation. However, the relation developed in Theorem 1 will provide a means of determining the false-alarm probability over the signal parameter domain from the detection probability.

The detection probability can be determined from the distribution of the LLR peak represented by (59) and the threshold in the form (54). A signal at $x$ with SNR $r = r(x)$ produces a peak log likelihood

$$\begin{aligned} \lambda_{\mathbf{m}} &= \langle \lambda_{\mathbf{m}} \rangle + \tilde{r}z \\ &= (r^2 + d)/2 + \tilde{r}z, \\ \tilde{r}^2 &= r^2 + d/2, \end{aligned} \tag{62}$$

where we have introduced a unit-variance variable $z$: the extra quadratic variance in the second term of (59), (which, as noted in Theorem 3, is half a chi-square variable of $d$ degrees of freedom), has been absorbed into the random term containing $z$. This approximation is legitimate in the asymptotic high-SNR domain under consideration ($r^2 >> d/2$), and we shall ignore the slight alteration from Gaussian in the compound distribution of $z$, since the added variance is small and the focus for the detection analysis is on the main body of the distribution.

We now infer from (62) that the probability of exceeding the local threshold $\lambda_r(x)$ in the presence of signal is

$$\begin{aligned} p_D(x) &= \text{prob}[\tilde{r}z \ge \lambda_r - (r^2 + d)/2] \\ &= \frac{1}{\sqrt{2\pi}} \int\limits_{z_D}^{\infty} e^{-z^2/2} dz, \\ z_D &= \tilde{r}^{-1}(\lambda_r - d/4) - \tilde{r}/2. \end{aligned} \tag{63}$$

As shown in Fig. 1, this probability is small for $d < r^2/2 < \lambda_r$, becomes appreciable for $r^2/2 \approx \lambda_r$ and closely approaches unity for $r^2/2 >> \lambda_r$. This has the effect of dividing the signal parameter space into regions of high and low probability separated by marginal regions of moderate probability, which we will soon show to contain most of the false alarms.

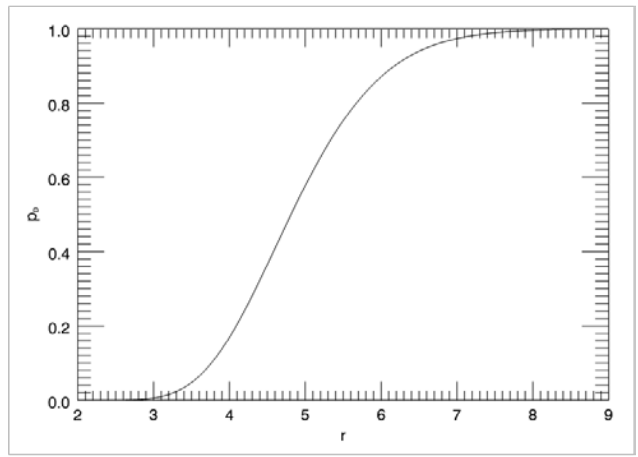

Figure 1. Probability of detection as a function of amplitude SNR $r(x)$; LLR threshold $\lambda_r$ (mid range value of 12.5) corresponds to 2-D simulation model to be introduced in Section 6.

A false alarm would appear as a spontaneous local maximum over threshold in $\lambda(x')$. For prediction, the mapping of the measurement, as a random quantity in its vector space, to the signal parameter domain through the generation of the LLR at each value of $x'$ creates a random field $\lambda(x')$. Predicting the occurrence distribution of local maxima over a multidimensional random field is a difficult problem that yields manageable analytic expressions only in the limit of asymptotically high values of a Gaussian field. Hasofer [15] presents a derivation of the asymptotic formula for the mean number of maxima above high levels in a Gaussian field, and Adler in his monograph, *The Geometry of Random Fields* [16], covers the analysis and provides many references to various approaches to the problem. But both of those authors treat only homogeneous fields ($r(x) = $ constant in our setting), and our primary interest is in a situation in which the SNR varies appreciably. Fortunately detection and false alarm quantities are connected through the differential relation of Theorem 1 in our optimized Bayes approach, and we can use that relation to obtain an analytic asymptotic expression for the false alarm density over the signal parameter domain even in the nonhomogeneous case.

We begin by manipulating the integrals on the right side of (27) to isolate the detection probability in a form suitable for our use. Under our assumption that prior probabilities and the geometric quantities (curvature) change only negligibly over the region $x \in X(x')$, we can evaluate them at $x = x'$, recognizing that $p_D(x'|x,\mu)$ is negligibly small outside that region. The reward function $R(x',x)$ is unity there [see (51)], and Equation (27) for the false alarm density can be written as

$$\nu_F(x') = \frac{-a(x')}{a_0 C(x')} \int\int_1^\infty \frac{1}{\mu} \frac{\partial}{\partial\mu} p_D(x'|x,\mu) d\mu \, dx. \qquad (64)$$

A change of variables from $\mu$ to $\lambda$,

$$\mu = \frac{a(x')V(x')}{a_0 C(x') r^d(x')} e^\lambda = e^{\lambda - \lambda_r}, \qquad (65)$$

connects to the amplitude log-likelihood threshold of (54):

$$\nu_F(x') = -\frac{r^d(x')}{V(x')} \int \int_{\lambda_r(x')}^\infty e^{-\lambda} \frac{\partial}{\partial\lambda} p_D(x'|x,\lambda) d\lambda \, dx. \qquad (66)$$

The detection quantity in the integral above is not the detection probability of (63), but it can be constructed from relations we have available. As identified after (23), it is the density over the signal parameter domain for detection at $x'$ of a signal at $x$ (over a log-likelihood threshold $\lambda$), and differentiation by $\lambda$ creates a density in log-likelihood amplitude as well. We represent this density as

$$-\frac{\partial}{\partial\lambda} p_D(x'|x,\lambda) = p_D(\lambda,x'|x) \qquad (67)$$

and decompose it into factors that we can evaluate:

$$p_D(\lambda',x'|x) = p_D(\lambda'|x',x) p_{x'}(x'|x). \qquad (68)$$

The prime on $\lambda$ indicates that it is the peak log-likelihood value at $x'$, obtained from $\lambda = \lambda(x)$ as indicated by (59):

$$\lambda' = \lambda + b, \qquad (69)$$

where the bias $b$,

$$b = \frac{1}{2}(x'-x) * r^2 Q(x'-x), \qquad (70)$$

is deterministic when $x'$ and $x$ are both conditions. Both $p_D(\lambda'|x',x)$, the probability density for observing a log-likelihood peak value of $\lambda'$ conditional on its occurring at a location $x'$ from a signal at $x$, and $p_{x'}(x'|x)$, the probability density that a log-likelihood peak value from a signal at $x$ occurs at $x'$, are rigorously Gaussian distributions (under the conditions of our analysis) and can be evaluated from their distribution parameters. The basic quantities needed are the covariances $\text{cov}(\lambda,x'|x)$ and $\text{cov}(x'|x)$, available from the random elements of (31), (44), and (58) (with the use of (36)) and directly from (60):

$$\text{cov}(\lambda,x'|x) = \frac{\nabla r*}{r} Q^{-1}$$
$$\text{cov}(x'|x) = \frac{1}{r^2} Q^{-1}. \qquad (71)$$

The gradient operator is represented as in (56) but with the differentiation occurring at $x$. The Gaussian distribution parameters can now be determined in the standard way:

$$\langle\lambda\rangle\big|_{x',x} = \langle\lambda\rangle\big|_x + \text{cov}(\lambda,x'|x)\big[\text{cov}(x'|x)\big]^{-1}(x'-x)$$
$$= \frac{1}{2} r^2(x) + \nabla*\left(\frac{1}{2} r^2(x)\right)(x'-x) \qquad (72)$$
$$\approx \frac{1}{2} r^2(x') = \frac{1}{2} r'^2$$

$$\text{var}(\lambda|x',x) = \text{var}(\lambda|x) -$$
$$\text{cov}(\lambda,x'|x)\big[\text{cov}(x'|x)\big]^{-1}\text{cov}(x',\lambda|x) \qquad (73)$$
$$= r^2(x)(1-\zeta) = \frac{1}{1+\gamma} r^2,$$

where the parameters

$$\gamma = \frac{\nabla r*}{r} M^{-1} \frac{\nabla r}{r}$$
$$\zeta = \frac{\nabla r*}{r} Q^{-1} \frac{\nabla r}{r} = \frac{\gamma}{1+\gamma} \qquad (74)$$

reflect the inhomogeneity of the field. Note that when the log-likelihood peak occurs at $x'$, the conditionally expected value of $\lambda$ at $x$ is its unconditionally expected value at $x'$ for a signal at $x'$. The conditional distribution parameters for $\lambda'$ are simply determined by incorporating the bias $b$, and we can construct the conditional density as

$$p_D(\lambda'|x',x) = \frac{\sqrt{1+\gamma}}{\sqrt{2\pi}\, r} \exp\left(-\frac{(1+\gamma)}{2r^2}\left[\lambda' - \frac{1}{2} r'^2 - b\right]^2\right) \qquad (75)$$

The integration indicated in the inner integral of (66) can now be performed (after the algebraic exercise of "completing the square" in $\lambda$) to yield the expression,

$$v_F\left(x'\right) = \frac{r'^d}{V}\exp\left(-\frac{1}{2}\zeta r'^2\right)\int G\left(z_r\left(x\right)\right)e^{-b}\,p_{x'}\left(x'\middle|x\right)dx\,, \quad (76)$$

where

$$z_r\left(x\right) = \frac{\left(1+\gamma\right)^{1/2}}{r'}\left(\lambda_r - \left[\frac{1}{2}\frac{\gamma-1}{\gamma+1}r'^2 + b\right]\right) \quad (77)$$

and

$$G\left(z_r\right) = \frac{1}{\sqrt{2\pi}}\int_{z_r}^{\infty}e^{-z^2/2}dz \quad (78)$$

is the cumulative Gaussian probability over the threshold $z_r$. Note that the $x$ dependence of $z_r$ arises from the bias term $b$ of (70). In evaluating the integral, we have substituted $r' = r\left(x'\right)$ for $r\left(x\right)$, an approximation that seems acceptable since the primary implication of a variable SNR has been incorporated through the parameters $\zeta$ and $\gamma$ (nonzero $\nabla r$). The distribution of the vector $x'$ is immediately available from its covariance (71) (continuing the substitution of $r'$ for $r$),

$$p_{x'}\left(x'\middle|x\right) = \frac{r'^d\left\|Q\right\|^{1/2}}{\left(2\pi\right)^{d/2}}e^{-b}\,, \quad (79)$$

and a change of variables from $x$ to the vector $z$ of independent, standard Gaussian components,

$$z = \left(r'^2 Q\right)^{1/2}\left(x'-x\right), \quad (80)$$

simplifies the expression for the false alarm density:

$$v_F\left(x'\right) = \frac{r'^d}{V}e^{-\zeta r'^2/2}\int G\left(z_r\left(b\right)\right)e^{-b}\frac{1}{\left(2\pi\right)^{d/2}}e^{-b}dz\,. \quad (81)$$

The integrand depends on $x$, and now $z$, only through the bias term $b$, which is evaluated as

$$b = \frac{1}{2}z^*z = \frac{1}{2}\left|z\right|^2\,, \quad (82)$$

and suggests a conversion to polar coordinates in the $d$ dimensions. With $u$ as the radial component, $u = \left|z\right|$, we can integrate over the $d$-1 angular components to produce $S_d$, the total angular measure (or the surface area of the unit $d$-dimensional sphere),

$$S_d = \frac{2\pi^{d/2}}{\Gamma\left(d/2\right)}\,, \quad (83)$$

where $\Gamma\left(n\right)$ is the gamma (factorial) function ($\Gamma\left(n+1\right) = n!$). (Note that $S_2 = 2\pi$ and $S_3 = 4\pi$ for familiar dimensions.)

*Theorem 4:* Our final expression for the false alarm density now assumes the form,

$$v_F\left(x'\right) = \frac{2\cdot 2^{-d/2}}{\Gamma\left(d/2\right)}\frac{r'^d}{V\left(x'\right)}e^{-\zeta r'^2/2}\int_0^\infty G\left(z_r\left(u\right)\right)e^{-u^2}u^{d-1}du\,, \quad (84)$$

where $z_r$ from (77) is now given by

$$z_r\left(u\right) = \frac{\left(1+\gamma\right)^{1/2}}{r'}\left(\lambda_r - \left[\frac{1}{2}\frac{\gamma-1}{\gamma+1}r'^2 + \frac{1}{2}u^2\right]\right)\,. \quad (85)$$

If the integration variable $u$ is replaced by the chi-square variable $u^2$, the gamma density for the chi-square distribution appears explicitly in (84) [14], and the integral (including the numerical coefficient) becomes recognizable as the expected value of $G\exp\left(-u^2/2\right)$.

We now have a form for the false alarm density that could hardly be obtained directly from the ambient distribution for the LLR. However, the form does not readily reveal qualitative characteristics of the density. Therefore, we make one additional approximation to obtain an expression for the false-alarm density that reveals the basic characteristics of the density over the signal parameter domain with more transparency but with less accuracy. In the next section we shall explore the distribution further with a 1-D example and compare it with simulation results.

We return to Equation (66) but now perform the integration over $x$ first. We assume that the SNR $r\left(x\right)$ as well as the prior probability and geometric properties are constant over the region $x\in X\left(x'\right)$ and neglect the correlations between $\lambda$ and $x'$ that result from a nonzero $\nabla r$ [see (71)]. Then the relation between signal location and detection location exhibited in (58) and (60) yields an interchange symmetry, $p_D\left(x\middle|x',\lambda\right) = p_D\left(x'\middle|x,\lambda\right)$, and the integration over $x$ in (66) becomes

$$\int p_D\left(x'\middle|x,\lambda\right)dx = \int p_D\left(x\middle|x',\lambda\right)dx\,. \quad (86)$$

But integration over the density variable leaves a quantity that is recognizable as the detection probability of (63). A change of amplitude variables (based on (63)) in the remaining integral,

$$z = \tilde{r}^{-1}\left(\lambda - d/4\right) - \tilde{r}/2\,, \quad (87)$$

leads to the form,

$$
\begin{aligned}
v_F(x') &= -\frac{r'^d e^{-d/4}}{V(x')} \int_{z_D}^{\infty} e^{-\bar{r}z - \bar{r}^2/2} \frac{\partial}{\partial z} \left\{ \frac{1}{\sqrt{2\pi}} \int_{z}^{\infty} e^{-z'^2/2} dz' \right\} dz \\
&= \frac{r'^d e^{-d/4}}{V(x')} \frac{1}{\sqrt{2\pi}} \int_{z_D}^{\infty} e^{-(z+\bar{r})^2/2} dz \\
&= \frac{r'^d e^{-d/4}}{V(x')} \frac{1}{\sqrt{2\pi}} \int_{z_F}^{\infty} e^{-z^2/2} dz \,,
\end{aligned}
\tag{88}
$$

where from (63)

$$
z_F = z_D + \bar{r} = \bar{r}^{-1}\left(\lambda_r - d/4\right) + \bar{r}/2 \,.
\tag{89}
$$

This expression, except for the factor $\exp(-d/4)$ and a somewhat compensatory $d$-dependence in the threshold and in $\bar{r}$ [see (62)], is the exact expression given in References [15] and [16] for the parameter-space density of maxima above a threshold $y_T = \lambda_r / r + r/2$ in a homogeneous standard Gaussian field of $d$ dimensions (the ambient field $y(x)$ of Equations (30) and (31)).

The qualitative behavior of the false alarm density with SNR in (88) is readily determined through an examination of the integration limit. The limit $z_F$ is minimum, and the density maximum, near $\bar{r}^2/2 = \lambda_r$, and the density declines rapidly for higher and lower SNR. Note that for both the expression above and the more general expression (84) the density at constant LLR threshold $\lambda_r(x')$ and constant SNR is inverse to the likelihood volume scale, showing that $V(x')$ is a measure of the intrinsic independent exposure to false alarms as discussed after Equation (53). But in fact the optimized threshold $\lambda_r$ also depends on the volume scale [see (54)] in a way that effectively cancels the explicit dependence on $V(x')$ to produce an approximate constant false alarm rate (CFAR) detector with regard to the volume scale dependence. This behavior, calculated from (88) and based on the 2-D signal model to be introduced in Section 6, is illustrated in Fig. 2. This figure also illustrates the difference between our threshold effect, based on minimization of a Bayes risk criterion, and that of the generalized likelihood ratio test (GLRT). Moustakides [7] presents conditions under which the GLRT is optimal (discrete estimation choices, constant Bayes costs and priors); in contrast, we show that it is distinctly nonoptimal in a situation in which the signal parameter dependence confers a strong $V(x')$ variation across the signal parameter domain. The GLRT approach lacks the integration over the continuous Bayes signal density that produces the volume dependence in the threshold, and it would result in a false alarm dependence following the dashed curve of Fig. 2.

The theoretical projections on false alarm behavior and detection estimation accuracy developed above are supported with simulations in the following section.

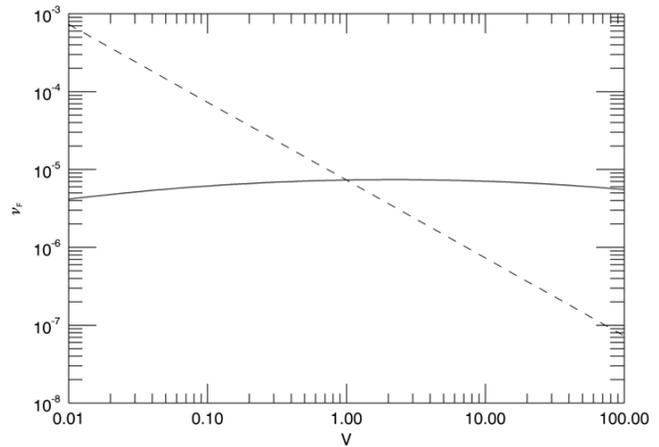

Figure 2. Relative false alarm density (88) as a function of the scale volume (53). Solid curve corresponds to optimal threshold $\lambda_r$ (54); dashed curve to threshold $\lambda_r$ with $V$ held constant ($V$=1).

## 6. SIMULATIONS

We apply a simple signal model and measurement setting to assess the accuracy of the theoretical expressions for false alarm densities ((84) and (88)) and state estimation accuracy (60) through Monte Carlo simulations. The model has two state dimensions, scalar position $\xi_s$ and signal width $\sigma$,

$$
x = \begin{pmatrix} \xi_s \\ \sigma \end{pmatrix},
\tag{90}
$$

and components of the signal vector $\mathbf{s}$ take the form,

$$
s_\xi(x) = A \exp\left(-\left(\xi - \xi_s\right)^2 / 2\sigma^2\right),
\tag{91}
$$

where $A$ is a known constant and the component index $\xi$ runs over integer values. In the simulations the noise content in the components $m_\xi$ of measurement $\mathbf{m}$ [see (1)] consists of independent draws of standard Gaussian variables, so the measurement covariance matrix $\mathbf{C}$ is the identity matrix. The measurement consists only of noise ($\mathbf{m} = \mathbf{n}$) in our simulations to test the theoretical false alarm expressions, and it includes a signal with specified state parameters in the simulations to test state estimation accuracy. Furthermore, for the false alarm assessment the model above will be reduced to two separate models of a single unknown parameter by specifying either $\xi_s$ or $\sigma$ as known—a reduction necessitated by computational constraints.

The LLR formulation of (29), the fundamental element of the analysis, applies a signal at a hypothesized state $x'$ to the measurement as

$$\lambda(x') = \mathbf{s}(x') \cdot \mathbf{m} - \frac{1}{2} r^2(x'), \qquad (92)$$

where the first term in (92), the matched filter output, is evaluated as

$$MFO(x') = \mathbf{s}(x') \cdot \mathbf{m} = \sum_{\xi} s_{\xi}(x') m_{\xi} \qquad (93)$$

and the SNR as

$$r^2(x') = \mathbf{s}(x') \cdot \mathbf{s}(x') = \sum_{\xi} \left[ s_{\xi}(x') \right]^2. \qquad (94)$$

Although computations for the theory-simulation comparison are performed as indicated in (93) and (94), the inner product for the SNR can be approximated through integration in a way that reveals the basic parameter dependence:

$$r^2(x') \approx A^2 \int \exp\left( -(\xi - \xi_{\mathbf{s}})^2 / \sigma^2 \right) d\xi = A^2 \sqrt{\pi}\sigma. \qquad (95)$$

This expression is actually very accurate except when the signal width parameter $\sigma$ is small: the fractional error is less than $10^{-3}$ when $\sigma > 1$.

Certain other quantities that depend on the form of the signal model (91) must be determined to perform the theory-simulation comparison. These include the LLR threshold and parameters in the theoretical expressions. The $Q$ matrix leads directly to several of these quantities, and for our signal and measurement models, it can be calculated from the relation in (55) as

$$Q = \frac{1}{r^2} \mathrm{cov}\left( \nabla \lambda \right) = \frac{1}{r^2} \begin{pmatrix} \dfrac{\partial \mathbf{s}}{\partial \xi_s} \cdot \dfrac{\partial \mathbf{s}}{\partial \xi_s} & \dfrac{\partial \mathbf{s}}{\partial \xi_s} \cdot \dfrac{\partial \mathbf{s}}{\partial \sigma} \\ \dfrac{\partial \mathbf{s}}{\partial \sigma} \cdot \dfrac{\partial \mathbf{s}}{\partial \xi_s} & \dfrac{\partial \mathbf{s}}{\partial \sigma} \cdot \dfrac{\partial \mathbf{s}}{\partial \sigma} \end{pmatrix}. \qquad (96)$$

The matrix expression on the right arises through the recognition from (92) that the random part of $\nabla \lambda$ is $(\nabla \mathbf{s}) \cdot \mathbf{n}$.

A short analytical exercise utilizing continuous integration as in (95) provides the $Q$ matrix in the approximate form,

$$Q \approx \begin{pmatrix} 1/2\sigma^2 & 0 \\ 0 & 3/4\sigma^2 \end{pmatrix}. \qquad (97)$$

Here, as throughout the simulation, the actual calculations were performed directly as vector inner products. The relative SNR variation with the state parameters is similarly approximated by

$$\frac{\nabla r}{r} = \frac{1}{r} \begin{pmatrix} \partial r / \partial \xi_s \\ \partial r / \partial \sigma \end{pmatrix} \approx \begin{pmatrix} 0 \\ 1/2\sigma \end{pmatrix} \qquad (98)$$

and the inhomogeneity parameters become, using the second equation of (74),

$$\zeta \approx 1/3 \\ \gamma \approx 1/2. \qquad (99)$$

For false alarm simulation an extremely large number of Monte Carlo trials is needed to produce a statistically significant number of false alarms above a high threshold. Therefore, for computational feasibility we reduced the 2-D state model to 1-D: first with $\sigma$ as the state variable and $\xi_s$ set to zero; then with $\xi_s$ as the state variable and $\sigma$ set to $\sigma_0 = 4$. Note that the 2-D model is simply adopted as a convenient source for the 1-D models. Results from the 1-D models have no direct bearing on false alarm expectations from the scenario with two unknown parameters. The full 2-D model is used to assess the estimation prediction of Theorem 3.

The threshold that is applied to the LLR values produced during the Monte Carlo trials of the false alarm simulation depends on what the state variable is and what its value is. With $\sigma$ as the sole state variable, the volume scale factor of (53), needed for the evaluation of $\lambda_r$, becomes

$$V(\sigma) = \sqrt{\frac{2\pi}{Q_{\sigma\sigma}}} \approx 2\sqrt{\frac{2\pi}{3}}\,\sigma \quad . \qquad (100)$$

While we have given approximate expressions in Equations (97) through (100), these quantities were rigorously calculated through the inner products of (96) and were used in (84) for comparison of the predicted density with the Monte Carlo results. With the parameter dependence of $r(x')$ and $V(x')$, the threshold $\lambda_r$ of (54) can be written as

$$\lambda_r(\sigma) = \lambda_0 + \log\left( r(\sigma)/r(\sigma_0) \right) - \log\left( V(\sigma)/V(\sigma_0) \right) \qquad (101)$$

and, using (95) and (100), in approximate form as

$$\lambda_r(\sigma) \approx \lambda_0 - \frac{1}{2}\log\left( \sigma/\sigma_0 \right), \qquad (102)$$

where the dependence on constants (including $A$), costs, and the *a priori* probabilities are lumped together in the constant $\lambda_0$. It is to $\lambda_0$ that we assign the various threshold values identified in the figures and table.

With $\sigma$ as the state variable, our Monte Carlo trials comprised random draws of standard Gaussian variates for the components of $N_{trial} = 5 \times 10^9$ measurement vectors **m** with the number of components sufficient to contain all computationally significant values of the signal form (91), and to each we applied the hypothesized signal (filter) as in (92) over a wide range of $\sigma$ values spaced by 0.2. The $\sigma$ values were searched, after interpolation to a finer grid (0.002 spacing), to find the one that maximized the LLR (92) for each Monte Carlo trial. These maximizing values, together with the associated maxima $\lambda(\sigma)$, were placed in an array over the finely spaced grid for comparison to a range of LLR thresholds $\lambda_r(\sigma)$ (101). The number of values above threshold in the finely spaced sigma bins were summed over a width $\Delta\sigma = 0.202$, centered about each point of the fine grid with 50 bins on either side, as a smoothing procedure. The number of values over this width ($N_\sigma$) was scaled to produce a simulation false alarm density as

$$\nu_F(\sigma)\big|_{\text{sim}} = \frac{1}{\Delta\sigma}\frac{N_\sigma}{N_{trial}} \ . \tag{103}$$

The corresponding theoretical density was computed from (84) with the threshold taken from (101) and the parameters $\zeta$ and $\gamma$ calculated from (74) for the characteristics of our model.

Results are presented in Figs. 3 to 5 for various LLR thresholds $\lambda_0$. While the densities were computed as functions of the state variable $\sigma$, they are plotted against the corresponding SNR amplitude $r$ (94) to substantiate our claim that false alarm production is maximized where the detection probability is marginal (where $r \approx \sqrt{2\lambda_0}$). To provide a clear confirmation of this behavior, we repeated the simulation for a different known signal amplitude $A$ (91) to produce Fig. 6. Compare Figs. 4 and 6 with different signal amplitude factors $A$ but identical LLR thresholds $\lambda_0$. False alarm production peaks at essentially the same SNR; it is clear from (95) that $\sigma$ adjusts to a different $A$ value to produce the same SNR peak location for the density; the density values change somewhat, but the peak occurs at essentially the same SNR. The predicted density agrees very well with the observed density in all cases; the maximum value appears less by just a small factor while the absolute density changes by about four orders of magnitude over the LLR threshold range of our calculations ($\lambda_0 = 5$ to 14).

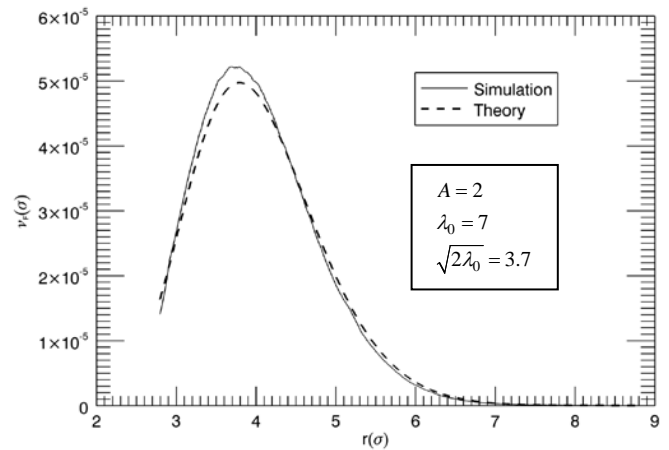

Figure 3. Comparison of false alarm densities in $\sigma$ for threshold $\lambda_0 = 7$. Simulation from (103); theory from (84).

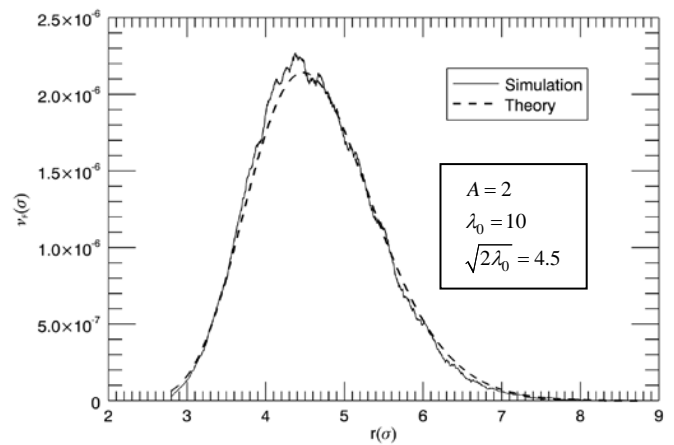

Figure 4. Comparison of false alarm densities in $\sigma$ for threshold $\lambda_0 = 10$. Simulation from (103); theory from (84).

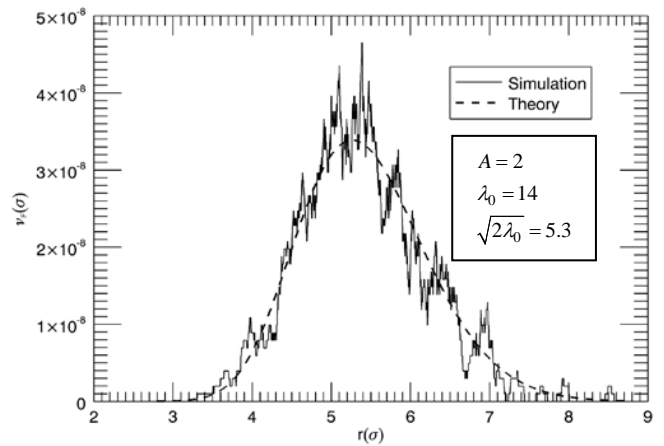

Figure 5. Comparison of false alarm densities in $\sigma$ for threshold $\lambda_0 = 14$. Simulation from (103); theory from (84).

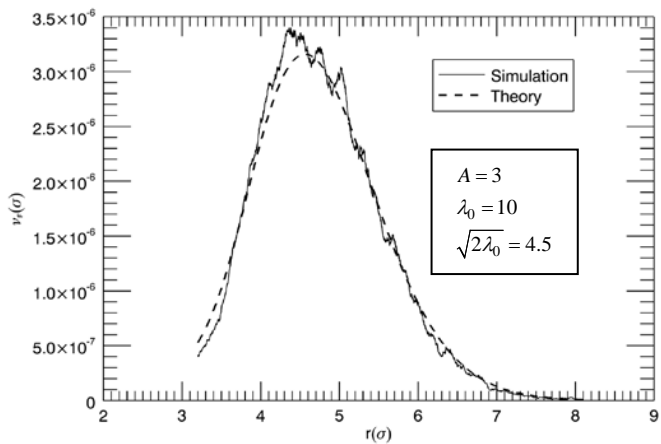

Figure 6. Comparison of false alarm densities in σ for threshold $\lambda_0 = 10$. Simulation from (103); theory from (84). Signal/filter amplitude factor changed from 2 to 3.

Integration of the densities over $\sigma$ provides a prediction of the total false alarm probability, and comparison of simulation and theoretical values at the various thresholds is presented in Table 1. The simulation values are actually obtained directly as the total count over the $\sigma$ domain of maximizing sigma values above threshold divided by the number of trials $N_{trial}$. The theoretical values are obtained through numerical integration of the density in (84) over $\sigma$. The agreement is clearly excellent despite the slight distortion noted above in the form of the densities

Table 1. Total false alarm probability – false alarm densities in $\nu_F(\sigma)$ integrated over σ domain. Parameter values: $A = 2$; $\sigma_0 = 4$.

| Threshold | Integrated FA Rate | |
|---|---|---|
| $\lambda_0$ | Expected (78) | Simulation |
| 5 | 7.05E-04 | 7.06E-04 |
| 6 | 2.81E-04 | 2.83E-04 |
| 7 | 1.08E-04 | 1.09E-04 |
| 8 | 4.08E-05 | 4.10E-05 |
| 9 | 1.52E-05 | 1.52E-05 |
| 10 | 5.63E-06 | 5.62E-06 |
| 11 | 2.08E-06 | 2.09E-06 |
| 12 | 7.67E-07 | 7.77E-07 |
| 13 | 2.83E-07 | 2.81E-07 |
| 14 | 1.04E-07 | 1.03E-07 |

The second 1-D false alarm comparison, with $\xi_s$ as the state variable, involves a homogeneous domain since the SNR is independent of $\xi_s$. Furthermore, if we limit our search over the domain to integer values of $\xi_s$, we can perform the search very efficiently with the use of the Fast Fourier Transform (FFT). The signal of (91) can be represented as $s_{\xi - \xi_s}$, symmetric in its subscript, and the matched filtering operation of (93) takes the form of a convolution:

$$MFO_{\xi_s} = \sum_\xi s_{\xi_s - \xi} m_\xi \ , \qquad (104)$$

The FFT of this expression yields the FFT of *MFO* as the product of the FFTs of *s* and *m*, and the inverse transform then produces the matched filter output of (104) for all hypothesized integer state locations $\xi_s$.

For this homogeneous case we generated measurement vectors **m** of high dimensionality (20,000 realizations, each with 10,000,000 $\xi$ samples) and processed each measurement vector as indicated by (104) using Fourier techniques to produce the matched filter output $\mathbf{s \cdot m}$ as a function of $\xi_s$. Local maxima of $\lambda$ (92) in $\xi_s$ were identified as potential false alarms, and a false alarm density was obtained by dividing the number of local maxima over threshold by the length of the domain (number of $\xi$ samples). Results were ensembled over the number of measurement vector realizations to provide additional statistical stability. This observed density, $\nu_F(\xi_s)|_{sim}$, is, of course, independent of $\xi_s$; because of homogeneity it applies to the entire $\xi_s$ domain.

To introduce an SNR dependence into the false alarm observation, we scale the signal model (91) over a range of amplitude factors $A$. Note that our original interpretation of $A$ remains: $A$ is not a state variable – we apply no operation to maximize $\lambda$ over $A$. The variation in $A$ simply represents a variation in the known signal amplitude. The same Monte Carlo trial realizations obtained for one $A$ value can be readily adapted to other values. The local maxima in the matched filter output occur at the same $\xi_s$ locations; the amplitudes simply scale with $A$. And the SNR term that is subtracted to form the LLR (92) scales as $A^2$. These scaling were used to obtain the plots of Figs. 7 and 8 for two different thresholds. Note that the abscissa on these plots is not $A$ but the SNR amplitude as determined from $A$ through (95). The applied threshold is unchanged as $A$ varies: $\lambda_r = \lambda_0$. The density of (88) is the appropriate theoretical comparison standard, and it provides the theoretical prediction in the plots. Again the agreement between theory and simulation is very good. While the general level of false alarm density changes by nearly two orders of magnitude between the two $\lambda_0$ threshold values of 10 and 14, at each threshold theory and simulation agree within about 5 to 10 percent. The theory-versus-simulation comparison could have been carried out equally well using the general formula (84) with the inhomogeneity parameters set to zero ($\gamma = \zeta = 0$); our calculations show that, although they involve slightly different approximations, the two theoretical expressions agree within one percent for these cases. The advantage of the expression in (88) for a homogeneous domain is that the algebraic form reveals the parametric dependence more clearly.

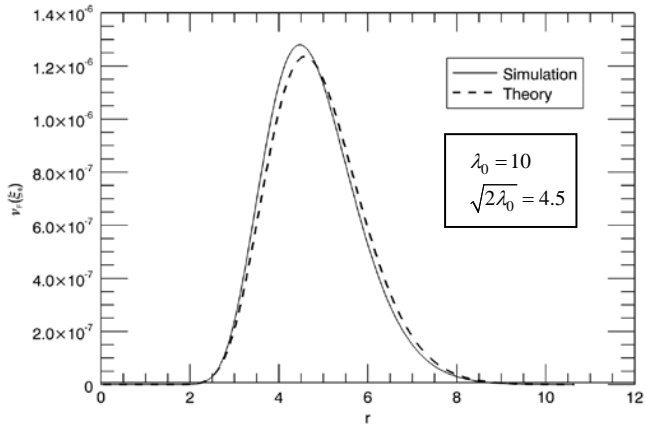

Figure 7. Homogeneous setting. Comparison of false alarm densities in $\xi_s$ over hypothesized SNR amplitudes. Theory from (88), threshold $\lambda_0 = 10$.

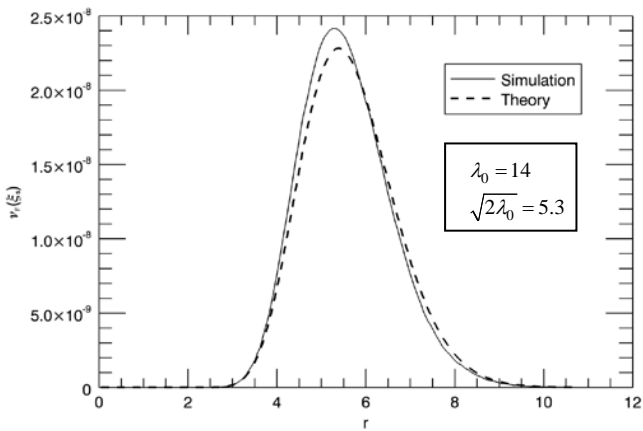

Figure 8. Homogeneous setting. Comparison of false alarm densities in $\xi_s$ over hypothesized SNR amplitudes. Theory from (88), threshold $\lambda_0 = 14$.

Our final simulation tests the accuracy prediction of Theorem 3 (60). The signal (91) was injected into a large number of measurements ( $\mathbf{m} = \mathbf{s} + \mathbf{n}$ ) with $\sigma = 4$ and with $\xi_\mathbf{s}$ selected uniformly between -0.5 and 0.5 over Monte Carlo trials (to avoid any unwanted bias in the digital processing over $\xi$ samples). The amplitude $A$ of the injected signal was varied to provide an SNR dependence in the accuracy test. The signal (91) was applied as a filter to the measurement in a 2-D search over $\sigma$ and $\xi_\mathbf{s}$ (over approximately $\pm 6$ standard deviations in each parameter about the signal injection location) to maximize $\lambda$ (92). Discrepancies between the state estimate and the injected state were tabulated and summarized in a sample covariance for comparison with the Cramer-Rao bound of Theorem 3. Comparisons of the diagonal entries, the variances for $\sigma$ and $\xi_s$, appear in Figs. 9 and 10. The approach of the simulation variances to the Cramer-Rao bound with increasing SNR is clear. The off-diagonal elements, which are zero theoretically, scattered about zero over the various SNR values of the simulation. They typically showed a sample correlation in the tenths of a percent, and the correlations were less than 1.5% for all SNR values.

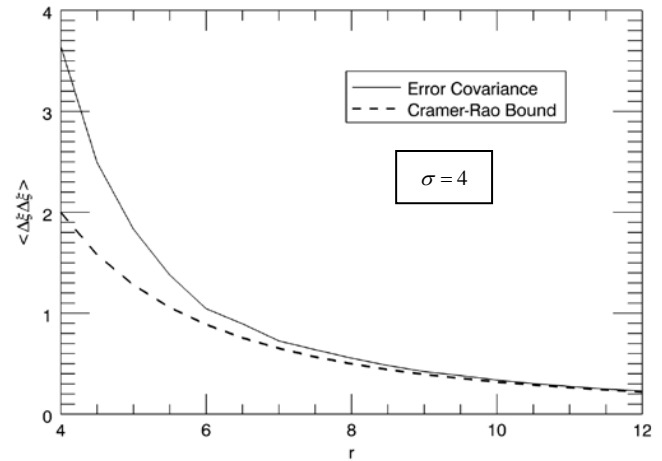

Figure 9. Comparison of covariance of $\xi_s$ estimate to Cramer-Rao bound as SNR amplitude varies.

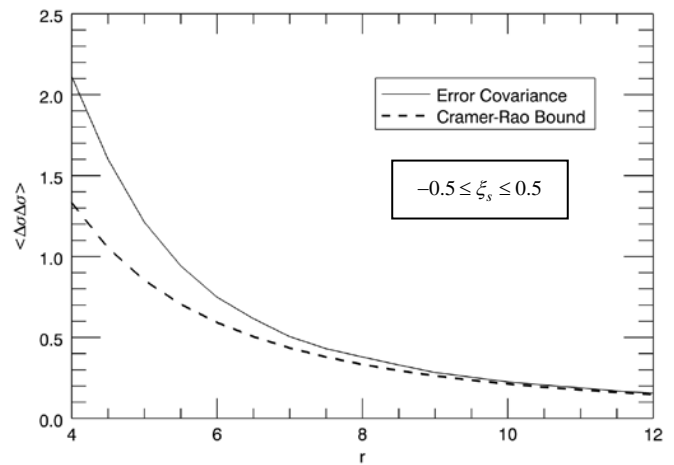

Figure 10. Comparison of covariance of $\sigma$ estimate to Cramer-Rao bound as SNR amplitude varies.

## 7. GLOBAL OPERATING CHARACTERISTIC

In a representative application, knowing the detection probability $p_D(x)$ as a function over the domain locations $x$ would typically be of significant value, but on the false alarm side, the probability of a false alarm anywhere in the parameter domain might be the quantity of most interest. It is possible, of course, to integrate the density $\nu_F(x)$ over the parameter domain, but the integration expression reveals little about the basic dependence of false alarm probability on threshold, and the integration itself is somewhat numerically awkward with the involvement of numbers extending over a very large dynamic range. However, the stationarity relation between detection and false alarm densities presented in Theorem 1 can be integrated over the search domain in $x$ to provide a global operating characteristic. In an environment in which the SNR varies significantly over the parameter domain $x$, this operating characteristic can be used to obtain a simple analytic expression for the global false alarm probability.

With our focus on global measures for detection and false alarms, we incorporate the $x$-dependence of Bayesian costs and rewards together with that of the signal prior probability $a(x)$ into a "pseudo" prior $a'(x)$. Constant factors can be absorbed into the threshold $\mu$ (which is varied to form the operating characteristic) so that $a'(x)$ is normalized to unity for the signal-present condition: $\int a'(x)dx = 1$. Using the approximations applied to obtain the expression (64) from Theorem 1 (27), the false alarm density corresponding to a threshold $\mu$ becomes

$$\nu_F(x',\mu) = -\int_\mu^\infty \frac{1}{\mu'} \frac{\partial}{\partial\mu'} \int a'(x') p_D(x'|x,\mu') dx\, d\mu' . \quad (105)$$

We first perform the integration over $x'$ to obtain the global false alarm probability $P_F(\mu)$,

$$P_F(\mu) = \int \nu_F(x',\mu) dx', \quad (106)$$

as

$$P_F(\mu) = -\int_\mu^\infty \frac{1}{\mu'} \frac{\partial}{\partial\mu'} \int a'(x) p_D(x,\mu') dx\, d\mu' , \quad (107)$$

recognizing in the process that $a'(x') \approx a'(x)$ where $p_D(x'|x,\mu')$ is of significant value [see the discussion preceding (64)]. The net *a posteriori* probability of a detection over the signal parameter domain is readily identified as

$$P_D(\mu) = \int a'(x) p_D(x,\mu) dx \quad (108)$$

to reveal the basic relation between false alarm and detection probabilities.

*Theorem 5*: The net probability of a false alarm can be determined from the net *a posteriori* probability of detection through the relation

$$P_F(\mu) = -\int_\mu^\infty \frac{1}{\mu'} dP_D(\mu') . \quad (109)$$

While this relation is obtained very simply from Theorem 1, it provides ready access to a predicted false alarm probability that would be very difficult to calculate directly. The relation is more readily appreciated with a change of independent variables from $\mu$ to $\lambda = \log \mu$ in $P_F$ and $P_D$. A change of integration variables $\lambda' = \log\left(\mu'/\mu\right)$ with thresholds related by $\lambda_T = \log \mu$ then leads directly to the expression

$$P_F(\lambda_T) = -e^{-\lambda_T} \int_0^\infty e^{-\lambda'} \frac{dP_D}{d\lambda} \left(\lambda_T + \lambda'\right) d\lambda' . \quad (110)$$

Restrict now to the case of strongly varying SNR, in which the parameter space is partitioned into islands of nearly certain detectability surrounded by regions of negligible detectability. In this case, the net detection probability, being crudely equal to the fraction of the domain over which detection is nearly certain (with the weighting factor $a'(x)$ incorporated), is actually a slowly varying function of the threshold. Then the derivative $dP_D / d\lambda$ in the integrand can be expanded to first order about $\lambda_T$ and the integration performed to yield the analytic expression

$$P_F(\lambda_T) \cong -e^{-\lambda_T} \frac{dP_D}{d\lambda} \left(\lambda_T + 1\right) . \quad (111)$$

This expression shows that the primary behavior of the net false alarm probability is a simple exponential dependence on the effective log-likelihood threshold $\lambda_T$.

## 8. CONCLUSIONS

When the general properties of the optimum Bayes detector/estimator in a parametric continuum of signals are applied to the specific case of additive Gaussian noise, several interesting new results emerge. First, in the asymptotic limit of high decision threshold, the signal state estimate, obtained from the location of a local maximum in the log-likelihood ratio, has an accuracy that approaches the Cramer-Rao bound. The decision threshold, depending on the signal parametric form, can vary across the parameter space, and in its application to the log-likelihood ratio, can divide that space into regions of nearly certain and nearly impossible detection, according to whether the squared SNR is greater than or less than twice the threshold value. Interestingly, the regions of marginal detectability separating the regions of high and low detectability contain the preponderance of false alarms. This last finding corresponds to a general integral relation between detection probability and false alarm density, a relation that can be interpreted as the operating characteristic when expressed as a function of arbitrary threshold. This integral relation is of especial value for determination of the false alarm density over the signal parameter domain. The direct determination of false alarm probability is a formidably difficult problem, but detection probabilities are more readily accessible and provide a reasonable path to false alarm prediction through the integral relation. Simulation results demonstrate a remarkably accurate prediction of false alarms. The form of the theoretical false alarm density over signal parameter space is well matched in the simulation results, and the agreement in total false alarm probability is almost exact over a range of thresholds that causes the probabilities to vary over orders of magnitude. Finally, the operating characteristic, extended globally over the entire signal parameter domain, leads to a simple expression for the net false alarm probability as a function of the net detection probability. In a signal parameter domain of strongly varying SNR, where direct estimation of the false alarm probability would be very difficult, the expression can provide a reasonable false alarm estimate from a crude estimate of where a signal is detectable.

The optimum Bayes detector/estimator, as developed here for a high-SNR, parametrically determined signal in additive Gaussian noise, exhibits features that distinguish it from other common detection/estimation approaches. Its relation to two of these – the CFAR detector and the GLRT – are mentioned briefly in the discussion centered on Fig. 2. The variation of the Bayes threshold with regard to the parameter-space volume of a likelihood ratio peak produces a false alarm dependence similar to that of a CFAR detector. Such a dependence is notably missing in a GLRT applied to a continuum signal parameter domain.

## ACKNOWLEDGMENT

The authors are grateful for the support of Arete Associates and for the generous help of Philip Pearle of Hamilton College.

**D. Michael Milder** earned his BS degree in physics at the California Institute of Technology, Pasadena CA USA, 1959, and his PhD in physics at Harvard University, Cambridge, MA USA, 1968. He worked at Tetra Tech from 1967 to 1972, R&D Associates from 1972-1975, and joined Arete Associates, Northridge, CA USA as a founding partner in 1976. His research papers have appeared in the Astrophysical Journal, the Journal of the Acoustical Society of America, the Journal of Fluid Mechanics, Radio Science, and Waves in Random and Complex Media. Dr. Milder passed away December 12, 2009.

**Robert G. Lindgren** earned his BS degree in chemical engineering at the University of Minnesota, Minneapolis, MN USA in 1965 and his PhD in chemical engineering with a minor in applied mathematics at the California Institute of Technology, Pasadena CA USA in 1970. He worked at R&D Associates from 1972-1982 and joined Arete Associates, Northridge, CA USA in 1982.

**Morris M. Berman** earned his BS degree in Mathematics at the California Institute of Technology, Pasadena CA USA, 1975, and his MS in Engineering and Applied Science at UCLA, Los Angeles, CA USA, 1978. He worked briefly at R&D Associates from 1975-1976, and has been with Arete Associates, Northridge, CA USA since its founding in 1976.